\documentclass[11pt, reqno]{amsart}

\usepackage[margin=3cm]{geometry}
\usepackage{amsmath,amssymb}
\usepackage{cite}
\usepackage{booktabs}
\usepackage{array}
\usepackage{tabularx}
\usepackage{float}
\usepackage{amsfonts}
\usepackage{enumerate,mathtools,bbm,bm,tcolorbox,mathrsfs}
\mathtoolsset{showonlyrefs}
\usepackage{hyperref}

\newtheorem{theorem}{Theorem}[section]
\newtheorem{proposition}[theorem]{Proposition}
\newtheorem{lemma}[theorem]{Lemma}

\theoremstyle{remark}

\allowdisplaybreaks
\numberwithin{equation}{section}

\newcommand{\R}{\mathbb{R}}
\newcommand{\C}{\mathbb{C}}

\newcommand{\Z}{\mathbb{Z}}

\newcommand{\SF}{\mathscr{F}}

\newcommand{\SN}{\mathscr{N}}

\newcommand{\ph}{\varphi}

\newcommand{\si}{\sigma}

\newcommand{\De}{\Delta}

\newcommand{\wt}{\widetilde}

\newcommand{\I}{\infty}

\newcommand{\II}{\mathbbm{1}}

\newcommand{\ls}{\lesssim}

\title[Alternative proof of endpoint Strichartz estimates]{An alternative proof of the Keel--Tao endpoint Strichartz estimates}
\author[S. Hadama]{Sonae Hadama}
\address{The University of Osaka, Japan}
\email{hadama.sonae.sci@osaka-u.ac.jp}
\thanks{The author was supported by JSPS KAKENHI Grant Number 26KJ0203.}

\date{}

\begin{document}
\begin{abstract}
	In this paper, we give an alternative proof of the abstract endpoint Strichartz estimates proved by Keel and Tao. Our method is based on the near-far decomposition technique, which was developed in a recent work by Toyoshima and the author.
\end{abstract}
\maketitle
\section{Introduction}
It is very well known that solutions to the free Schr\"odinger equation satisfy \textit{Strichartz estimates}
\begin{equation}\label{eq:Schrodinger}
	\|e^{it\De} \ph\|_{L^p_t(\R;L^q_x(\R^d))} \ls \|\ph\|_{L^2_x(\R^d)}
\end{equation}
for $p,q\in [2,\I]$ such that
\[
\frac2p+\frac{d}q = \frac{d}2, \qquad (d,p,q)\ne (2,2,\I).
\]
Strichartz \cite{Strichartz 1977} first proved \eqref{eq:Schrodinger} for $p=q=(2d+4)/d$.
Later, Ginibre and Velo \cite{Ginibre Velo 1985, Ginibre Velo 1995} broadened the range of $(p,q)$ except the endpoint case $p=2$ in an abstract setting.
Finally, in \cite{Keel Tao 1998}, Keel and Tao resolved the endpoint case $p=2$ also in an abstract setting. More precisely, Keel and Tao established the following theorem.
\begin{theorem}[Keel--Tao]\label{thm:endpoint}
	Let $I\subset \R$ be a (possibly unbounded) interval.
	Let $(X,\mu)$ be a $\sigma$-finite measure space and let $H$ be a separable complex Hilbert space.
	Assume that $U(t):H\to L^2(X)$ is bounded and that $t\mapsto U(t)\varphi$ is continuous for every $\varphi\in H$.
	Suppose that $\si>1$ and 
		\begin{align}
			&\sup_{t\in I}
			\|U(t)\|_{H \to L^2(X)} \le B,\label{eq:energy}\\
			&\|U(t)U(\tau)^*\|_{L^1(X)\to L^\I(X)}
			\le C_0|t-\tau|^{-\sigma}.\label{eq:dispersion}
		\end{align}
	Then, the following estimates hold:
	\begin{align}
		&\|U(t)\varphi\|_{L^2(I;L^{\frac{2\sigma}{\sigma-1},2}(X))}
		\le C \|\varphi\|_H,\label{eq:homogeneous} \\
		&\bigg\|\int_{\tau<t}U(t)U(\tau)^*F(\tau)d\tau\bigg\|_{L^2(I;L^{\frac{2\sigma}{\sigma-1},2}(X))}
		\le C^2 \|F\|_{L^2(I;L^{\frac{2\sigma}{\sigma+1},2}(X))},\label{eq:retarded}
	\end{align}
	where $L^{q,r}$ is the standard Lorentz space.
	Moreover, we may take $C=C_\sigma C_0^{\frac1{2\sigma}}B^{1-\frac1\sigma}$, where $C_\sigma>0$ only depends on $\sigma$.
\end{theorem}

Keel and Tao \cite{Keel Tao 1998} proved the endpoint estimates using a dyadic decomposition in time, with the key summation step justified either by an explicit atomic decomposition or by bilinear real interpolation.
Vilela \cite{Viela 2001} gave an alternative proof of the endpoint Strichartz estimate for the free Schrödinger equation with radial initial data.
Dodson \cite{Dodson 2019} presented a new proof for the free Schr\"odinger equation based on level-set and Whitney decompositions.

In this paper, we will give an alternative proof of Theorem \ref{thm:endpoint}.
Our proof adapts a \textit{near-far decomposition} developed in \cite{Hadama Toyoshima 2026+}.
In \cite{Hadama Toyoshima 2026+}, this method was used to prove the orthonormal Strichartz estimates.
We will show that this method is also useful to prove (standard) Strichartz estimates.

An advantage of the present proof is that we can keep the spatial exponents fixed at the endpoint when we sum the estimates.
Keel and Tao divide the time integral into dyadic pieces.
Since the endpoint bound is the same at every time scale, they use estimates at nearby spatial exponents to sum these pieces.
In our proof, for each pair of dyadic levels $k,k'$ of the squared weight $f^2$, it is enough to split the time integral into two parts.
We use the near part $0<t-\tau\leq\ell(k,k')$ and the far part $t-\tau>\ell(k,k')$. With this choice of cutoff, the bounds for both parts decay geometrically in $|k-k'|$.

\subsection{Notation}
To shorten the notation, we often write $L^p_t L^q_x$ for $L^p_t (I;L^q(X))$.
For $f:I \times X \to \C$, define
\[
(fU\ph)(t,x):= f(t,x)(U(t)\ph)(x) \qquad \text{for } \ph\in H,
\]
and its formal adjoint by
\[
(fU)^*F:= \int_I U(t)^* (\overline{f(t)}F(t)) dt \qquad \text{for } F:I\times X \to \C.
\]
For $f,g:I \times X \to \C$, define
\[
(fUU^*g)_<: F(t)\mapsto f(t)U(t)\int_{\tau<t} U(\tau)^*(g(\tau)F(\tau))d\tau.
\]

\section{An alternative proof of Theorem \ref{thm:endpoint}}
We first prove the following key proposition.
\begin{proposition}\label{prop:weighted}
Under the same assumptions as in Theorem \ref{thm:endpoint},
\begin{equation}\label{eq:weighted}
\begin{aligned}
 &\|(fUU^*g)_<\|_{L^2(I\times X)\to L^2(I\times X)} \le C_\sigma C_0^{\frac1\sigma}B^{2-\frac2\sigma}\|f\|_{L^\infty(I;L^{2\sigma,\infty}(X))} \|g\|_{L^\infty(I;L^{2\sigma,\infty}(X))}.
\end{aligned}
\end{equation}
\end{proposition}

\begin{proof}
\noindent \textbf{Step 1: Near-far decomposition.}
First, we assume that $f=g\ge 0$ and
\begin{equation}\label{eq:levels}
	V:=f^2 =\sum_{k\in\mathbb Z}2^k\II_{E_k},
\end{equation}
where $(E_k)_{k\in\Z}$ is disjoint and $E_k=\varnothing$ except for finitely many $k\in \Z$.
We write $E_k(t)=\{x:(t,x)\in E_k\}$.
Define
\begin{equation}\label{eq:time-scale}
 \ell(k,k'):= C_0^{\frac1\sigma}B^{-\frac2\sigma} \|V\|_{L^\infty_t L^{\sigma,\infty}_x} 2^{-\frac{k+k'}{2}}2^{\frac{|k-k'|}{4}}.
\end{equation}
Then, we have
\begin{align}
	&f(t)U(t)\int_{\tau<t}U(\tau)^* (f(\tau) u(\tau))d\tau\\
	&\qquad = \sum_{k,k'\in \Z} \II_{E_k}(t)f(t)U(t)\int_{\tau<t}U(\tau)^* (f(\tau) \II_{E_{k'}}(\tau)u(\tau))d\tau \\
	&\qquad = \sum_{k,k'\in \Z} \II_{E_k}(t)f(t)U(t)\int_{t-\tau> \ell(k,k')}U(\tau)^* (f(\tau) \II_{E_{k'}}(\tau)u(\tau))d\tau \\
	&\qquad \qquad + \sum_{k,k'\in \Z} \II_{E_k}(t)f(t)U(t)\int_{0<t-\tau\le \ell(k,k')}U(\tau)^* (f(\tau) \II_{E_{k'}}(\tau)u(\tau))d\tau \\
	&\qquad =: \mathscr{F} + \mathscr{N}.
\end{align}

\noindent \textbf{Step 2: Estimate of the far part $\SF$.}
By the definition of the weak norm, for almost every $t$,
\begin{equation}\label{eq:level-measure}
 \mu(E_k(t))
 \le2^{-\sigma k}\|V\|_{L^\infty_t L^{\sigma,\infty}_x}^\sigma.
\end{equation}
By the dispersion assumption \eqref{eq:dispersion},
\begin{equation}\label{eq:spatial-far}
 \begin{aligned}
 \big\|\II_{E_k(t)}f(t)U(t)U(\tau)^*
        f(\tau)\II_{E_{k'}(\tau)}\big\|_{L^2_x\to L^2_x} 
        &\le C_0|t-\tau|^{-\sigma}2^{\frac{k+k'}{2}}
             \mu(E_k(t))^{\frac12}\mu(E_{k'}(\tau))^{\frac12}\\
 &\le C_0|t-\tau|^{-\sigma}2^{-\frac{(\sigma-1)(k+k')}{2}}
             \|V\|_{L^\infty_t L^{\sigma,\infty}_x}^\sigma.
 \end{aligned}
\end{equation}
Since $\sigma>1$,
\begin{equation}\label{eq:time-integral}
 \int_{\ell(k,k')}^\infty \frac1{\tau^{\sigma}}d\tau
 =\frac{\ell(k,k')^{1-\sigma}}{\sigma-1}.
\end{equation}
Thus, Young's inequality in time and \eqref{eq:time-scale} give
\begin{align}
&\bigg\|\II_{E_k(t)}f(t)U(t)\int_{t-\tau>\ell(k,k')}U(\tau)^*\bigl(f(\tau)\II_{E_{k'}(\tau)}u(\tau)\bigr)d\tau\bigg\|_{L^2_{t,x}}\\
 &\qquad\le C_\sigma C_0 2^{-\frac{(\sigma-1)(k+k')}{2}}
 \ell(k,k')^{1-\sigma}
 \|V\|_{L^\infty_t L^{\sigma,\infty}_x}^\sigma
 \|\II_{E_{k'}}u\|_{L^2_{t,x}}\\
 &\qquad= C_\sigma C_0^{\frac1\sigma} B^{2-\frac2\sigma}
 \|V\|_{L^\infty_t L^{\sigma,\infty}_x} 2^{-\frac{(\sigma-1)|k-k'|}4}\|\II_{E_{k'}}u\|_{L^2_{t,x}}.
 \end{align}
Therefore, we obtain
\begin{align}
	 \|\mathscr{F}\|_{L^2_{t,x}}
	 &\le \bigg\|\II_{E_k(t)}f(t)U(t)\int_{t-\tau>\ell(k,k')}U(\tau)^*\bigl(f(\tau)\II_{E_{k'}(\tau)}u(\tau)\bigr)d\tau\bigg\|_{\ell_k^2 \ell_{k'}^1 L^2_{t,x}} \\
	 &\le C_\sigma C_0^{\frac1\sigma} B^{2-\frac2\sigma}
	 \|V\|_{L^\infty_t L^{\sigma,\infty}_x}
	 \Big\|2^{-\frac{(\sigma-1)|k-k'|}4}\|\II_{E_{k'}}u\|_{L^2_{t,x}}\Big\|_{\ell_k^2 \ell_{k'}^1} \\
	 &\le C_\sigma C_0^{\frac1\sigma} B^{2-\frac2\sigma}
	 \|V\|_{L^\infty_tL^{\sigma,\infty}_x}\|u\|_{L^2_{t,x}}.
\end{align}

\noindent \textbf{Step 3: Estimate of the near part $\SN$.}
We have
\begin{align}
	\|\mathscr{N}\|_{L^2_{t,x}}
	&\le \bigg\|\II_{E_k(t)}f(t)U(t)\int_{0<t-\tau\le \ell(k,k')}U(\tau)^*\bigl(f(\tau)\II_{E_{k'}(\tau)}u(\tau)\bigr)d\tau\bigg\|_{\ell_k^2 \ell_{k'}^1 L^2_{t,x}}.
\end{align}
When $k\ge k'$, the energy bound \eqref{eq:energy} gives
\begin{align}
&\bigg\|\II_{E_k(t)}f(t)U(t)\int_{0<t-\tau \le \ell(k,k')}
 U(\tau)^*\bigl(f(\tau)\II_{E_{k'}(\tau)}u(\tau)\bigr)d\tau\bigg\|_{L^2_{t,x}} \\
&\qquad \le \bigg\|\int_{0<t-\tau \le \ell(k,k')}
\|f(t)U(t)U(\tau)^*(f(\tau)\II_{E_{k'}(\tau)}u(\tau))\|_{L^2_x} d\tau\bigg\|_{L^2_t} \\ 
&\qquad \le \ell(k,k')^{\frac12} \|f(t)U(t)U(\tau)^*(f(\tau)\II_{E_{k'}(\tau)}u(\tau))\|_{L^2_t L^2_\tau L^2_x} \\
&\qquad \le B \|fU\|_{H \to L^2_{t,x}} 2^{\frac{k'}2} \ell(k,k')^{\frac12}\|\II_{E_{k'}}u\|_{L^2_{t,x}}\\
&\qquad = \|fU\|_{H\to L^2_{t,x}} C_0^{\frac1{2\sigma}} B^{1-\frac1\sigma} \|V\|_{L^\infty_t L^{\sigma,\infty}_x}^{\frac12} 2^{-\frac{|k-k'|}8} \|\II_{E_{k'}}u\|_{L^2_{t,x}}.
\end{align}
When $k< k'$, the energy bound \eqref{eq:energy} also gives
\begin{align}
	&\bigg\|\II_{E_k(t)}f(t)U(t)\int_{0<t-\tau \le \ell(k,k')}
	U(\tau)^*\bigl(f(\tau)\II_{E_{k'}(\tau)}u(\tau)\bigr)d\tau\bigg\|_{L^2_{t,x}} \\
	&\qquad \le B 2^{\frac{k}2} \bigg\|\int_{0<t-\tau \le \ell(k,k')}
	U(\tau)^*(f(\tau)\II_{E_{k'}(\tau)}u(\tau))d\tau\bigg\|_{L^2_t H} \\
	&\qquad \le B \|(fU)^*\|_{L^2_{t,x} \to H} 2^{\frac{k}2}\|\II_{E_{k'}} u\|_{L^2_t L^2_\tau([t-\ell(k,k'),t];L^2_x)} \\
	&\qquad \le B \|fU\|_{H \to L^2_{t,x}} 2^{\frac{k}2} \ell(k,k')^{\frac12} \|\II_{E_{k'}} u\|_{L^2_{t,x}} \\
	&\qquad = \|fU\|_{H\to L^2_{t,x}} C_0^{\frac1{2\sigma}} B^{1-\frac1\sigma} \|V\|_{L^\infty_t L^{\sigma,\infty}_x}^{\frac12} 2^{-\frac{|k-k'|}8} \|\II_{E_{k'}}u\|_{L^2_{t,x}}.
\end{align}
Therefore, we obtain
\begin{align}
	\|\mathscr{N}\|_{L^2_{t,x}}
	&\le \|fU\|_{H\to L^2_{t,x}} C_\sigma C_0^{\frac1{2\sigma}} B^{1-\frac1\sigma} \|V\|_{L^\infty_t L^{\sigma,\infty}_x}^{\frac12}
	     \Big\|2^{-\frac{|k-k'|}8} \|\II_{E_{k'}}u\|_{L^2_{t,x}} \Big\|_{\ell_k^2 \ell_{k'}^1} \\
	&\le \|fU\|_{H\to L^2_{t,x}} C_\sigma C_0^{\frac1{2\sigma}} B^{1-\frac1\sigma} \|V\|_{L^\infty_t L^{\sigma,\infty}_x}^{\frac12}
	\|u\|_{L^2_{t,x}}.
\end{align}

\noindent \textbf{Step 4: Conclusion.}
By the above arguments, 
\begin{align}
	&\|(fUU^*f)_< u\|_{L^2_{t,x}}
	\le \|\mathscr{F}\|_{L^2_{t,x}} + \|\mathscr{N}\|_{L^2_{t,x}} \\
	&\qquad \le C_\sigma C_0^{\frac1\sigma} B^{2-\frac2\sigma} \|V\|_{L^\infty_t L^{\sigma,\infty}_x} \|u\|_{L^2_{t,x}}
	           + \|fU\|_{H\to L^2_{t,x}} C_\sigma C_0^{\frac1{2\sigma}} B^{1-\frac1\sigma} \|V\|_{L^\infty_t L^{\sigma,\infty}_x}^{\frac12}
	           \|u\|_{L^2_{t,x}}.
\end{align}
Hence,
\begin{align}
	A&:=\|(fUU^*f)_<\|_{L^2_{t,x}\to L^2_{t,x}}\\
	&\le C_\sigma C_0^{\frac1\sigma} B^{2-\frac2\sigma} \|V\|_{L^\infty_t L^{\sigma,\infty}_x} 
	+ \|fU\|_{H\to L^2_{t,x}} C_\sigma C_0^{\frac1{2\sigma}} B^{1-\frac1\sigma} \|V\|_{L^\infty_t L^{\sigma,\infty}_x}^{\frac12}.
\end{align}
Since
\begin{align}\label{eq:half weighted}
	\|fU\|_{H\to L^2_{t,x}}^2
	= \|fUU^* f\|_{L^2_{t,x}\to L^2_{t,x}}
	\le 2 \|(fUU^* f)_<\|_{L^2_{t,x} \to L^2_{t,x}},
\end{align}
we have
\begin{align}
	A&\le C_\sigma C_0^{\frac1\sigma} B^{2-\frac2\sigma} \|V\|_{L^\infty_t L^{\sigma,\infty}_x}
	+ \sqrt{A} C_\sigma C_0^{\frac1{2\sigma}} B^{1-\frac1\sigma} \|V\|_{L^\infty_t L^{\sigma,\infty}_x}^{\frac12}\\
	&\le C_\sigma C_0^{\frac1\sigma} B^{2-\frac2\sigma} \|V\|_{L^\infty_t L^{\sigma,\infty}_x}
	+ \frac12 A.
\end{align}
Since $E_k=\varnothing$ except for finitely many $k\in \Z$, we can see $A<\infty$ and we obtain \eqref{eq:weighted} when $V$ is a dyadic step function.
Note that we can remove the assumption $\#\{k:E_k\neq \varnothing\}<\I$.
Finally, we prove \eqref{eq:weighted} for general $f, g\ge 0$. 
Let
\[
\|f\|_{L^\infty_t L^{2\sigma,\infty}_x} = \|g\|_{L^\infty_t L^{2\sigma,\infty}_x} = 1.
\]
Let $h:=\max(f,g)$ and choose a dyadic step function $\wt{h}^2=\wt{V}$ such that
\[
 \frac12\wt{h} \le h \le \wt{h}.
\]
Hence,
\begin{align}
	\|(fUU^*g)_<\|_{L^2_{t,x}\to L^2_{t,x}}
	&\le \bigg\|\frac{f}{\wt{h}}\bigg\|_{L^\infty_{t,x}} \bigg\|\frac{g}{\wt{h}}\bigg\|_{L^\infty_{t,x}}
	\Big\|(\wt{h}UU^* \wt{h})_< \Big\|_{L^2_{t,x}\to L^2_{t,x}}
	\le C_\sigma C_0^{\frac1\sigma}B^{2-\frac2\sigma}\|\wt{V}\|_{L^\infty_t L^{\sigma,\infty}_x}.
\end{align}
By the definition of the weak norm,
\[
\|\wt{V}\|_{L^\infty_t L^{\sigma,\infty}_x} \le C_\sigma.
\]
Therefore, we obtain \eqref{eq:weighted}.
\end{proof}

We now prove Theorem \ref{thm:endpoint}.
\begin{proof}[Proof of Theorem \ref{thm:endpoint}]
By \eqref{eq:half weighted} and Proposition \ref{prop:weighted},
\begin{align}
	\|U(t)\varphi\|_{L^2_t L^{\frac{2\sigma}{\sigma-1},2}_x}
	&\le C \sup_{\|f\|_{L^\infty_t L^{2\sigma,\infty}_x} \le 1} \|f(t)U(t)\varphi\|_{L^2_t L^2_x} \\
	&\le C\|\varphi\|_H \sup_{\|f\|_{L^\infty_t L^{2\sigma,\infty}_x} \le 1}
	     \|(fUU^*f)_<\|_{L^2_{t,x} \to L^2_{t,x}}^{\frac12}\le C_\sigma C_0^{\frac1{2\sigma}}B^{1-\frac1\sigma} \|\varphi\|_H,
\end{align}
which gives \eqref{eq:homogeneous}.
For any decomposition $F=F_0 F_1$,
\begin{align}
	&\bigg\|U(t)\int_{\tau < t} U(\tau)^* F(\tau) d\tau\bigg\|_{L^2_t L^{\frac{2\sigma}{\sigma-1},2}_x} \\
	&\qquad \le C \sup_{\|f\|_{L^\infty_t L^{2\sigma,\infty}_x} \le 1}
	             \bigg\|f(t)U(t)\int_{\tau < t} U(\tau)^* F(\tau) d\tau\bigg\|_{L^2_t L^2_x} \\
	&\qquad \le C \sup_{\|f\|_{L^\infty_t L^{2\sigma,\infty}_x} \le 1}
	             \|(fUU^*F_0)_< \|_{L^2_{t,x} \to L^2_{t,x}} \|F_1\|_{L^2_{t,x}} \\
	&\qquad \le C_\sigma C_0^{\frac1\sigma}B^{2-\frac2\sigma} \|F_0\|_{L^\infty_t L^{2\sigma,\infty}_x} \|F_1\|_{L^2_{t,x}}.
\end{align}
By Lemma \ref{lem:factor}, we can choose $F_0, F_1$ such that 
\[
\|F_0\|_{L^\infty_t L^{2\sigma,\infty}_x} \|F_1\|_{L^2_{t,x}} \le C_\si \|F\|_{L^2_t L_x^{\frac{2\sigma}{\sigma+1},2}}.
\]
Hence, we conclude that \eqref{eq:retarded} holds.
\end{proof}

\appendix

\section{Factorization of mixed Lorentz spaces}

The author believes that the following lemma is known, but has not been able to find a suitable reference. We therefore include a proof for completeness.

\begin{lemma}\label{lem:factor}
	Let $I\subset \R$ be a (possibly unbounded) interval and let $(X,\mu)$ be a $\si$-finite measure space.
	Let
	\[
	1\le p,q,r,s,p_i,q_i,r_i,s_i\le\infty,
	\qquad i=0,1,
	\]
	satisfy
	\[
		\frac1p=\frac1{p_0}+\frac1{p_1},
		\qquad \frac1q=\frac1{q_0}+\frac1{q_1},
		\qquad \frac1r=\frac1{r_0}+\frac1{r_1},
		\qquad \frac1s=\frac1{s_0}+\frac1{s_1}.
	\]
	Assume that $q_i=\I$ if $p_i=\infty$. Similarly, we also assume that $s_i=\I$ if $r_i = \I$.
	Then, for every $F \in L^{p,q}_t L^{r,s}_x$, there exist $F_0, F_1$ such that $F=F_0F_1$ and
	\[
	\|F_0\|_{L^{p_0,q_0}(I;L^{r_0,s_0}_x(X))}
	\|F_1\|_{L^{p_1,q_1}(I;L^{r_1,s_1}_x(X))}
	\ls \|F\|_{L^{p,q}(I;L^{r,s}(X))},
	\]
	where the implicit constant depends only on the exponents.
\end{lemma}

\begin{proof}
	Applying \cite[Proposition 2.1]{Lorist Nieraeth 2022} in $x$ for every fixed $t\in I$, we can decompose $F(t)=G_0(t) G_1(t)$ such that
	\[
	\|G_0(t)\|_{L^{r_0,s_0}(X)} \|G_1(t)\|_{L^{r_1,s_1}(X)} \le C \|F(t)\|_{L^{r,s}(X)}.
	\]
	Applying \cite[Proposition 2.1]{Lorist Nieraeth 2022} in $t$ again, we have $g_0$ and $g_1$ such that
	$g_0(t)g_1(t)=\|F(t)\|_{L^{r,s}(X)}$ and 
	\[
	\|g_0\|_{L^{p_0,q_0}(I)}
	\|g_1\|_{L^{p_1,q_1}(I)}
	\le C \|F\|_{L^{p,q}(I;L^{r,s}(X))}.
	\]
	Define
	\[
	F_0(t,x)=\frac{g_0(t)G_0(t,x)}{\|G_0(t)\|_{L^{r_0,s_0}(X)}},
	\qquad
	F_1(t,x)=\frac{\|G_0(t)\|_{L^{r_0,s_0}(X)}G_1(t,x)}{g_0(t)}.
	\]
	Then $F=F_0F_1$ and
	\begin{align*}
	&\|F_0\|_{L^{p_0,q_0}(I;L^{r_0,s_0}(X))}\|F_1\|_{L^{p_1,q_1}(I;L^{r_1,s_1}(X))} \\
	&\qquad \le \|g_0\|_{L^{p_0,q_0}(I)} \bigg\|\frac{\|G_0(t)\|_{L^{r_0,s_0}(X)}\|G_1(t)\|_{L^{r_1,s_1}(X)}}{g_0(t)}\bigg\|_{L^{p_1,q_1}(I)} \\
	&\qquad \le C \|g_0\|_{L^{p_0,q_0}(I)} \|g_1\|_{L^{p_1,q_1}(I)} \le C \|F\|_{L^{p,q}_t(I;L^{r,s}(X))},
	\end{align*}
	which completes the proof.
\end{proof}

\section*{Declaration of AI use}
ChatGPT Sol 6.1 and Astra 6 were used to assist in polishing proofs and in identifying typographical errors. All content was finally reviewed and checked by the author, who takes full responsibility for the final manuscript.


\begin{thebibliography}{9}

\bibitem{Dodson 2019}
B. Dodson,
\textit{Defocusing Nonlinear Schr\"odinger Equations},
Cambridge Tracts in Mathematics, vol. \textbf{217}, Cambridge University Press, Cambridge, 2019.

\bibitem{Ginibre Velo 1985}
J. Ginibre and G. Velo,
\textit{The global Cauchy problem for the non linear Schr\"odinger equation revisited},
Ann. Inst. H. Poincar\'e Anal. Non Lin\'eaire \textbf{2} (1985), no. 4, 309--327.

\bibitem{Ginibre Velo 1995}
J. Ginibre and G. Velo,
\textit{Generalized Strichartz inequalities for the wave equation},
J. Funct. Anal. \textbf{133} (1995), no.~1, 50--68.

\bibitem{Hadama Toyoshima 2026+}
S. Hadama and S. Toyoshima,
\textit{Endpoint orthonormal Strichartz estimates},
arXiv preprint (2026), arXiv:2608.21823.

\bibitem{Keel Tao 1998}
M. Keel and T. Tao,
\textit{Endpoint Strichartz estimates},
Amer. J. Math. \textbf{120} (1998), no. 5, 955--980.

\bibitem{Lorist Nieraeth 2022}
E. Lorist and Z. Nieraeth,
\textit{Sparse domination implies vector-valued sparse domination},
Math. Z. \textbf{301} (2022), no. 1, 1107--1141.

\bibitem{Strichartz 1977}
R. S. Strichartz,
\textit{Restrictions of Fourier transforms to quadratic surfaces and decay}
of solutions of wave equations,
Duke Math. J. \textbf{44} (1977), no. 3, 705--714.

\bibitem{Viela 2001}
M. C. Vilela,
\textit{Regularity of solutions to the free Schrödinger equation with radial initial data},
Illinois J. Math. \textbf{45} (2001), no. 2, 361--370.

\end{thebibliography}
\end{document}